\documentclass[12pt,twoside]{amsart}
\usepackage{amsmath,amsthm,amsfonts,amssymb,amscd,epsfig}

\def\ca{{\mathcal A}}
\def\cb{{\mathcal B}}
\def\cc{{\mathcal C}}

\def\cg{{\mathcal G}}

\def\cl{{\mathcal L}}

\def\ct{{\mathcal T}}

\def\h{\hspace{0.02in}}

\def\Br{\mathbb R}
\def\Bz{\mathbb Z}
\def\Bn{\mathbb N}

\newtheorem{thm}{Theorem}
\newtheorem{thm*}{Theorem}

\newtheorem{lem}{Lemma}
\newtheorem*{lem*}{Lemma}

\theoremstyle{definition}

\newtheorem*{defn*}{Definition}
\newtheorem{example}{Example}

\newtheorem*{rem*}{Remark}

\newtheorem*{cor*}{Corollary}
\newtheorem{prop*}{Proposition}

\theoremstyle{definition}

\newtheorem{exmp*}{Example}

\newtheorem{exmps*}{Examples}

\theoremstyle{remark}

\newtheorem*{rems*}{Remarks}

\newtheorem*{ack*}{Acknowledgment}

\begin{document}

\title
[On the definition of relative pressure]{A comment on the definition
  of relative pressure}

\author{Karl Petersen}
\address{Department of Mathematics,
CB 3250, Phillips Hall,
         University of North Carolina,
Chapel Hill, NC 27599 USA}
\email{petersen@math.unc.edu}

\author{Sujin Shin}
\address{
Department of Mathematics, Korea Advanced Institute of Science
and Technology, Daejon, 305-701, South Korea}
\email{sjs@math.kaist.ac.kr}

\begin{abstract}
We show that two natural definitions of the relative pressure function
for a locally constant potential function and a factor map from a shift
of finite type coincide almost everywhere with respect to every
invariant measure. With a suitable extension of one of the
definitions, the same holds true for any continuous potential function.
\end{abstract}
\maketitle

The introduction to the paper \cite{PQS} included, for factor maps
between subshifts and the identically zero potential, a
``finite-range'' definition of the relative pressure function that is
different from the standard one \cite{LW,Wal}, which involves complete
bisequences. While this variation in the definition had no bearing on
the results of that paper, it does seem useful to clarify the extent
to which the two definitions differ; in particular, in some situations
one definition may be easier to use than the other. We show in this
note that for a factor map $\pi : X \to Y$, where $X$ is a shift of
finite type and $Y$ is a subshift, the two relative pressure functions
can be different, but they coincide almost everywhere with respect to
every invariant measure on $Y$.  Therefore, for each ergodic invariant
measure $\nu$ on $Y$, the two definitions lead to the same value of
the maximal possible relative entropy $h_\mu (X|Y)$ of any invariant
measure $\mu$ on $X$ over $Y$. More generally, for any locally
constant potential function (one that depends on only finitely many
coordinates), the analogously defined two relative pressure functions
coincide almost everywhere with respect to every invariant measure on
$Y$.  Finally, we show that this statement continues to hold for an
arbitrary continuous potential function with a suitably generalized 
definition of the finite-range relative pressure function.

Several useful ideas for the study of the relative entropy function
were developed in \cite{Shi} in connection with questions about the
existence of compensation functions, and we adopt and extend them for
our purposes here.

If $(X, S)$ is a topological dynamical system, then $M (X)$  
will denote the set of all $S$-invariant Borel probability
measures on $X$ with its weak* topology. Given $x \in X$, let
\begin{equation*}
\mu_x = \lim\limits_{n \rightarrow \infty} \frac{1}{n} 
\sum\limits^{n - 1}_{i = 0} \delta_{S^i x} \in M (X)
\end{equation*}
if it exists, in which case we call $x$ a {\em generic point}. 
Denote by $\cg (X)$ the set of all generic points in $X$. 
If $X$ is a shift space, then for each $n \geq 1$, $\cb_n (X)$ denotes
the set of $n$-blocks in $X$, and $\cb (X) = \cup_n \cb_n (X)$.
Given $b_1 \cdots b_n \in \cb_n (X)$, $n \geq 1$, 
we define 
\begin{equation*}
_{r}[b_1 \cdots b_n]_{r+n-1} = \{x \in X : x_r=b_1, \dots
,x_{r+n-1}=b_n\},
\end{equation*}
and we abbreviate $_0[b_1 \cdots b_n]_{n-1}$ by $[b_1 \cdots b_n]$.
We denote the shift transformation by $\sigma$ and the usual metric
for a subshift by $\rho$.

Let $S : X \rightarrow X$ and
$T : Y \rightarrow Y$ be continuous maps of compact 
metrizable spaces and $\pi : X \rightarrow Y$ a factor map.
i.e., a continuous surjection with $\pi \circ S = T \circ \pi$.
For a given compact subset $K$ of $X$, for $n \geq 1$ and 
$\delta > 0$, denote by $\Delta_{n, \delta} (K)$ the set of 
$(n, \delta)$-separated sets of $X$ contained in $K$. 
Let $f \in C(X)$. Fix $\delta > 0$ and $n \geq 1$.
For each $y \in Y$, let
\begin{equation*}
P_n (\pi, f, \delta) (y) = \sup \bigg\{ \sum\limits_{x \in E}
\exp \Big( \sum_{i = 0}^{n - 1} f (S^i x) \Big) \Big\vert E \in
\Delta_{n, \delta} \big( \pi^{-1} \{ y \} \big) \bigg\}.
\end{equation*} 
Define $P (\pi, f) : Y \rightarrow \Br$ by
\begin{equation*}
P(\pi, f) (y) = \lim\limits_{\delta \rightarrow 0}
\limsup\limits_{n \rightarrow \infty} \frac{1}{n} \ln
P_n (\pi, f, \delta) (y).
\end{equation*}
The function $P (\pi, f)$ is called the 
{\em relative pressure function} associated with $f$. It is 
Borel measurable and $T$-invariant. 
For $\nu \in M (Y)$, let $M (\nu)=\pi^{-1}(\nu)$ denote the
set of measures in $M (X)$ that project to $\nu$ under $\pi$.
Given $f \in C(X)$, the function 
$P (\pi, f) : Y \rightarrow \Br$ satisfies the 
{\em relative variational principle} \cite{LW}: 
for each $\nu \in M(Y)$,
\begin{equation*}
\int\!P(\pi, f) \h d \nu = \sup\limits \bigg\{ h(\mu) 
+ \int\!f d \mu \bigg\vert \mu \in M(\nu) \bigg\} - h(\nu).
\end{equation*}
In particular, for a fixed $\nu \in M(Y)$, 
\begin{equation*}
\sup \{h_{\mu}(X|Y): \mu \in M(\nu)\} = \sup \{h(\mu)- h(\nu): \mu \in
M(\nu)\} = \int_Y P(\pi ,0)\,  d \nu .
\end{equation*}

Hereinafter, let $X$ be a shift of finite type and $Y$ a subshift, on
finite alphabets $\ca (X)$ and $\ca (Y)$, respectively.  
Let $\pi : X \rightarrow Y$ be a factor map, so that $Y$ is a sofic
system. We treat the 2-sided case, the 1-sided case being very
similar. 
Let $f \in \cc (X)$ be a locally constant function, i.e. one that
depends on only finitely many coordinates 
$x_{-m} \cdots x_m$.
By passing to a higher block representation if necessary, we may
assume that $\pi$ is 
 represented
by a one-block map from $\cb_1 (X)$ to $\cb_1 (Y)$, which we denote
again by $\pi$, and that $f$ is a function of the two coordinates
$x_0x_1$.

 For $B = b_1 \cdots b_n \in \cb_n (X)$, 
$\pi (B)$ means the $n$-block $\pi (b_1) \cdots \pi (b_n)$ of $Y$;
given $v \in \cb_n (Y)$, $\pi^{-1} (v)$ 
denotes the set of $n$-blocks 
of $X$ that project to $v$ by the block map $\pi$. 
Given $y \in Y$, for each $n \geq 1$,
let $D_n (y)$ consist of one point from each nonempty set 
$\pi^{-1} (y) \cap [x_0 x_1 \cdots x_{n - 1}]$.
The potential function $f$ determines 
a block map $F : \cb_2 (X) \rightarrow \Br$ by 
$F (b_0b_1) = \exp (f (x))$ for any $x \in [b_0b_1]$. 
For a block $B = b_1 \cdots b_n \in \cb (X)$, put 
\begin{equation*}
s_f (B) = F (b_1b_2)F(b_2b_3) \cdots F (b_{n-1}b_n),
\end{equation*}
and for a block $w \in \cb (Y)$, put $S_f (w) = \sum_B s_f (B)$ 
where the sum is taken over all $B \in \cb (X)$
that are mapped to $w$ by $\pi$. Then for each $y \in Y$,
\begin{equation}
\begin{split}
\label{eq_P(pi,f)}
P(\pi,f) (y) & = \limsup_{n \rightarrow \infty} \frac{1}{n} \ln 
\bigg[ \sum_{x \in D_n (y)} 
\exp \Big( \sum^{n-1}_{i=0} f (\sigma^i x) \Big) \bigg] \\
& = \limsup_{n \rightarrow \infty} \frac{1}{n} \ln
\bigg[ \sum_{x \in D_n (y)} s_f (x_0 x_1 \cdots x_{n - 1}) \bigg]
\end{split}
\end{equation} 
(see \cite[Theorem 4.6]{Wal}). 
In particular, for all $y \in Y$,
\begin{equation*}
P (\pi, 0) (y) = \limsup\limits_{n \rightarrow \infty} \frac{1}{n} 
\ln \big\vert D_n (y) \big\vert
\end{equation*}
(with $f \equiv 0$). Define another Borel-measurable function 
$\Phi_f : Y \rightarrow \Br$ by
\begin{equation*}
\Phi_f (y) = \limsup_{n \rightarrow \infty} \frac{1}{n} 
\ln S_f (y_0 y_1 \cdots y_{n - 1}) \hspace{0.3in} \text{for } y \in Y.
\end{equation*}
It can be shown that $\Phi_f (y) \leq \Phi_f (\sigma y)$ 
for all $y \in Y$.
Also $P (\pi, f) (y) \leq \Phi_f (y)$ for all $y \in Y$. 
Given $f \in C (X)$, one may have $P (\pi, f) (y) < \Phi_f (y)$ for
some  $y \in Y$, as seen in the following example (which in
\cite{Shi1} and \cite{Shi} was shown to be a factor map for which
there exists no saturated compensation function).

\begin{example}
\label{ex1}
Let $X$, $Y$ be the subshifts of finite type
determined by allowing the transitions marked on 
Figure \ref{fig:ex1} and the one-block factor code 
$\pi : X \rightarrow Y$ map $1$ to $1$, and $2, 3, 4, 5$ to $2$. 

\begin{figure}[h]
\centerline{\epsfig{figure=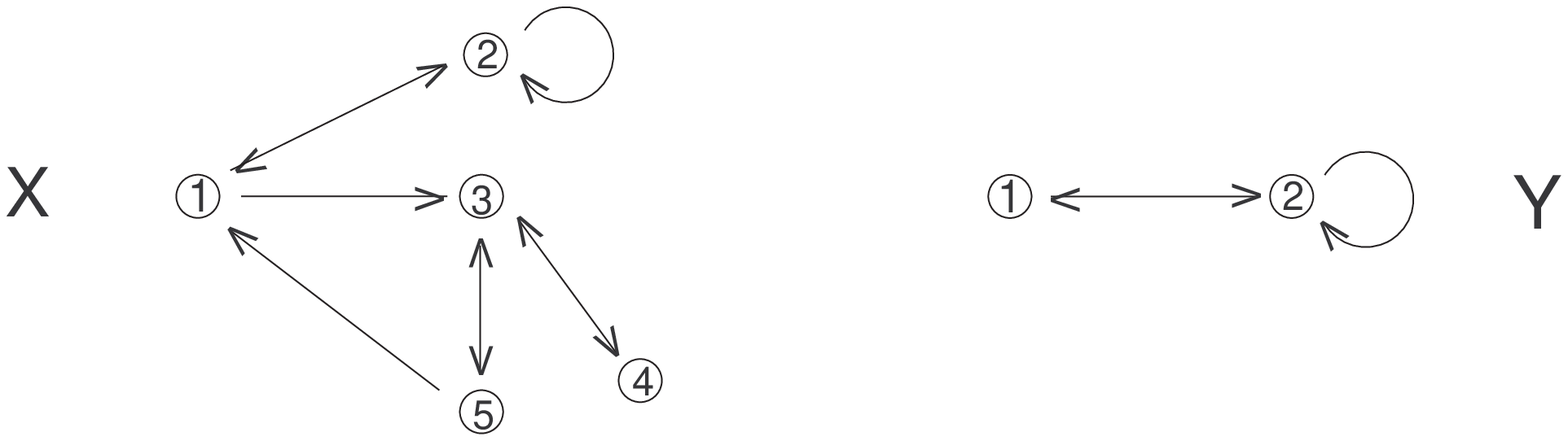, height=1.3in}}
\caption{}
\label{fig:ex1}
\end{figure}

For each $k \geq 1$, let
$a_k = 2^k + 1$, and define $y \in Y$ by
\begin{equation*}
y = \cdots 2 2 2 . 1 2^{a_1} 1 2^{a_2} 1 2^{a_3} 1 \cdots
\end{equation*}
(so that $y_i = 2$ for all $i < 0$).
Whenever $\pi x = y$, then $x_{[0, \infty)} = y_{[0, \infty)} 
= 1 2^{a_1} 1 2^{a_2} 1 2^{a_3} 1 \cdots$,
since each $a_k$ is odd and so $\pi^{-1} (1 2^{a_k} 1) 
= \{ 1 2^{a_k} 1 \}$. Thus $\vert D_n (y) \vert = 1$ for all 
$n \geq 1$ which implies that $P (\pi, 0) (y) = 0$. 
Meanwhile, fix $k \geq 1$ and set $n_k = 2^{k + 1} + 2 k - 2$. Then
\begin{equation*}
\big\vert \pi^{-1} [y_0 y_1 \cdots y_{n_k - 1}] \big\vert
= \big\vert \pi^{-1} [1 2^{a_1} 1 2^{a_2} 1 \cdots 1 2^{a_k}]
\big\vert = \big\vert \pi^{-1} [1 2^{a_k}] \big\vert 
= 2^{2^{k - 1}} + 1.
\end{equation*}
Consider $f \equiv 0 \in C (X)$. Then
\begin{equation*}
\begin{split}
\Phi_f (y) & \geq \limsup\limits_{k \rightarrow \infty} \frac{1}{n_k} 
\ln \big\vert \pi^{-1} [y_0 y_1 \cdots y_{n_k - 1}] \big\vert \\
& = \lim\limits_{k \rightarrow \infty} 
\frac{\ln \big( 2^{2^{k - 1}} + 1 \big)}{2^{k + 1} + 2 k - 2} 
= \frac{\ln 2}{4} > 0 = P (\pi, 0) (y). 
\end{split}
\end{equation*}
\end{example}

Our goal is to prove the following.

\begin{thm}
\label{thm-main}
Let $X$ be an irreducible shift of finite type, $Y$ a subshift, and 
$\pi : X \rightarrow Y$ a factor map. 
Let $f \in C (X)$ be a function which depends on
just $2$ coordinates, $x_0x_1$, of each point $x \in X$. 
Then for each $\nu \in M (Y)$, we have $P(\pi,f)(y) = \Phi_f(y)$
a.e. $d\nu (y)$, and hence
\begin{equation*}
\int\!P (\pi, f) \h d \nu = \int\!\Phi_f \h d \nu.
\end{equation*}
\end{thm}


\begin{rem*}
We may assume that $f \geq 0$. For if $f$ takes some negative values,
choose a constant $c$ such that $f+c \geq 0$ and use the equations
$P(\pi,f+c) = P(\pi,f)+c, \Phi_{f+c}=\Phi_f +c$. Thus we have $1 \leq
F \leq M$ for some constant $M$.
\end{rem*}

For notational convenience, set $s = s_f$,
$S = S_f$ and $\Phi = \Phi_f$, and let 
\begin{equation*}
\ct (y) = P (\pi, f) (y) \hspace{0.3in} \textnormal{for } y \in Y.
\end{equation*}
The map $\mathcal{R} : M (X) \rightarrow \Br^+$ 
defined by $\mathcal{R} (\mu) = h (\mu) - h (\pi \mu)$ 
is upper semicontinuous in the $\textnormal{weak}^*$»õ
topology \cite[Lemma 2.2]{Wal}. Using this one can easily prove the
following \cite{Shi}. 

\begin{lem}
\label{lem-cl_usc}
For any $f \in C(X)$ the 
affine map $\cl : M(Y) \rightarrow \Br$ given by
$\cl (\nu) = \int\!\ct \h d \nu$ is upper semicontinuous.  
\end{lem}

For $q \geq 1$, let $P_q (Y) = \{ y \in Y \vert \sigma^q (y) = y \}$ 
and $P (Y) = \bigcup_{q \geq 1} P_q (Y)$. 

\begin{lem}
\label{lem-P(Y)}
Let $y \in P_q (Y)$, $q \geq 1$. Then
\begin{equation}
\label{eq-nq}
\Phi (y) = \lim\limits_{n \rightarrow \infty} \frac{1}{nq} 
\ln S (y_0 y_1 \cdots y_{n q - 1}).
\end{equation}
Also $\Phi (y) = \ct (y)$.
\end{lem}

\begin{proof}
For a block $w \in \cb (Y)$ with $w w \in \cb (Y)$, 
we have $S_f (w^{n + m}) \leq M \cdot S_f (w^n) S_f (w^m)$. Thus 
$(1 / n) \ln S_f (w^n)$ converges as $n \rightarrow \infty$
(see \cite[p. 240]{Pet}),
and hence
\begin{equation*}
\Phi (y) = \lim_{n \rightarrow \infty} \frac{1}{nq} 
\ln S (y_0 \cdots y_{n q - 1}).
\end{equation*}

To show that $\Phi (y) = \ct (y)$, set 
$w = y_0 y_1 \cdots y_{q - 1} \in \cb_q (Y)$, so that
\begin{equation*}
y = \cdots w . w w \cdots 
= \cdots y_{q - 1} . y_0 y_1 \cdots y_{q - 1} y_0 y_1 \cdots.
\end{equation*}
Let $\vert \pi^{-1} (w) \vert = l \geq 1$ and 
define an $l \times l$, $0$-$1$ matrix $A = (A_{uv})$ by
$A_{u v} = 1$ if and only if $u v \in \cb_{2q} (X)$, where 
$u, v \in \pi^{-1} (w)$, that is,
$A$ is the transition matrix between blocks in
the inverse image of the repeating word $w$ that forms $y$. 
Note that some blocks in $\pi^{-1}(w)$ may be preceded or followed
only by allowable $q$-blocks in $X$ that do not map to $w$. So we let
$B$ be the reduction of $A$ obtained by excluding all the zero columns
and rows together with their corresponding rows and columns. Then 
by (\ref{eq-nq}), 
\begin{equation*}
\Phi (y) = \lim_{n \rightarrow \infty} \frac{1}{nq} 
\ln \bigg[ \sum_{u_i \in \pi^{-1} (w)}
s (u_1 \cdots u_n) A_{u_1 u_2} \cdots A_{u_{n - 1} u_n} \bigg].
\end{equation*}
Since $B$ is essential, we have 
\begin{equation*}
\begin{split}
\sum_{u_i \in \pi^{-1} (w)} s (u_1 \cdots 
& u_{n + 1} u_{n + 2}) A_{u_1 u_2} \cdots A_{u_{n + 1} u_{n + 2}} \\
& \leq \sum_{u_i \in \pi^{-1} (w)} \big[ s (u_1 \cdots u_n) 
B_{u_1 u_2} \cdots B_{u_{n - 1} u_n} \big] l^2 M^{2 q} \\
& \leq l^2 M^{2 q} \sum_{x \in D_{nq} (y)} 
s (x_0 x_1 \cdots x_{nq - 1}) \\ 
& \leq l^2 M^{2 q} \sum_{u_i \in \pi^{-1} (w)} 
s (u_1 \cdots u_n) A_{u_1 u_2} \cdots A_{u_{n - 1} u_n}.
\end{split}
\end{equation*}
Now we take logarithms, divide, and take 
the limit on $n$ to get $\Phi (y) = \ct (y)$.
\end{proof}

For the proof of the following result, we refer to \cite{Shi}.

\begin{lem}
\label{lem-Omega}
Let $y \in \cg(Y)$ and let $y^{(s)} \in P_{l_s} (Y)$, $l_s \geq s$, 
for each $s \geq 1$. If there is $N \geq 1$ such that
$y^{(s)}_{[0, l_s - N]} = y_{[0, l_s - N]}$ for all $s$ large enough,
then $\mu_{y^{(s)}} \rightarrow \mu_y$ as $s \rightarrow \infty$.  
\end{lem}

Let $\ca$ now denote the alphabet of $X$. 
Fix $y \in Y$. Given $b, c \in \ca$ and $n \geq 1$, let
\begin{equation*}
\Gamma^n_y (b, c) = \sum_u S (b u c) ,
\end{equation*}

where the sum is taken over all $u$'s in $\cb_{n - 1} (X)$
such that $\pi (b u c) = y_0 y_1 \cdots y_n$.
(If $\Gamma^n_y (b, c) \geq 1$, then $\pi b = y_0$ and 
$\pi c = y_n$.) Then
\begin{equation*}
\sum_{b, c \in \ca} \Gamma^n_y (b, c) = S (y_0 y_1 \cdots y_n).
\end{equation*}
It is not difficult to check the following.

\begin{lem}
\label{lem-Gamma}
Let $y \in P_q (Y)$, $q \geq 1$, and $b \in \ca$. Then for $k \geq 1$,  
\begin{equation*}
\big[ \Gamma^q_y (b, b) \big]^k \leq \Gamma^{qk}_y (b, b).
\end{equation*}
\end{lem}

\begin{lem}
\label{lem-mu_Phi}
Let $y \in \cg(Y)$. Then there is a sequence 
$\{ y^{(s)} \}^\infty_{s = 1} \subset P (Y)$ such that
$\mu_{y^{(s)}} \rightarrow \mu_y$ as $s \rightarrow \infty$ and 
$\Phi (y) \leq \liminf_{s \rightarrow \infty} \Phi (y^{(s)})$. 
\end{lem}

\begin{proof}
Observe first that there exist a symbol, say $a$, of $\ca(Y)$ and
a strictly increasing sequence 
$\{ m_k \}^\infty_{k = 0} \subset \Bz^+$ such that
$y_{m_k} = a$ for all $k \geq 0$ and 
\begin{equation}
\label{eq-m_k}
\Phi (y) = \lim_{m_k \rightarrow \infty} \frac{1}{m_k} 
\ln S (y_0 \cdots y_{m_k}).
\end{equation} 
Let $y^* = \sigma^{m_0} y \in Y$ and for each 
$k \geq 0$ put $n_k = m_k - m_0 \geq 0$. Then $y^*_i = a$ 
if $i = n_k$ for some $k \geq 0$ ($n_0 = 0$). 
Let $\cc = \pi^{-1} (a)$. 
For $k \geq 1$, choose $b_k, c_k \in \cc$ so that
\begin{equation*}
\Gamma^{n_k}_{y^*} (b_k, c_k) 
= \max_{b, c \in \cc} \Gamma^{n_k}_{y^*} (b, c)
\end{equation*} 
(so $\pi (b_k) = y^*_0$ and $\pi (c_k) = y^*_{n_k}$).
Since $1 \leq F \leq M$, it follows that
\begin{equation*}
\frac{1}{M^2} S (y^*_0 \cdots y^*_{n_k}) \leq 
\Gamma_{y^*}^{n_k} (b_k, c_k) \leq S (y^*_0 \cdots y^*_{n_k}).
\end{equation*}
Thus by (\ref{eq-m_k}),
\begin{equation}
\label{eq-Phi=lim} 
\begin{split}
\Phi (y) & \leq \liminf_{k \rightarrow \infty} \frac{1}{m_k} 
\ln \big[ M^{m_0} \cdot S (y_{m_0} \cdots y_{m_k}) \big] \\
& = \liminf_{k \rightarrow \infty} \frac{1}{n_k} 
\ln S (y^*_0 \cdots y^*_{n_k}) = \liminf_{k \rightarrow \infty} 
\frac{1}{n_k} \ln \big[ \Gamma^{n_k}_{y^*} (b_k, c_k) \big]. 
\end{split}
\end{equation}

Notice that there exist $b_*, c_* \in \cc$ such that
$b_k = b_*$ and $c_k = c_*$ for infinitely many $k$'s, say 
$k_s$'s, where $k_s \nearrow \infty$ as $s \rightarrow \infty$.
Since $X$ is irreducible, there is $w \in \cb_d (X)$ for some 
$d \geq m_0$ such that $c_* w b_* \in \cb_{d + 2} (X)$. 
Let $a_1 \cdots a_d = \pi (w) \in \cb_d (Y)$.
Fix $s \geq 1$ and put $l_s = n_{k_s} + 1 + d$.
Define $y^{(s)} \in P_{l_s} (Y)$ by
\begin{equation*}
y^{(s)} = \cdots a_d \h . \h y^*_0 y^*_1 \cdots y^*_{n_{k_s}} 
a_1 \cdots a_d \h y^*_0 y^*_1 \cdots
\end{equation*}
($y^*_0 \dots y^*_{n_{k_s}}$ is the image under $\pi$ of a word 
$w^* =w_0^* \cdots w^*_{n_{k_s}}$ in $X$ with 
$w_0^* =b_k=b_*$ and 
$w_{n_{k_s}}^*=c_k=c_*$, 
so that $\cdots w.w^*ww^*w\cdots$ is a legitimate point
$x^{(s)} \in X$, and $y^{(s)} = \pi (x^{(s)})$).
Since $y^{(s)}_{[0, l_s - d]} = y^*_{[0, l_s - d]}$ and 
$l_s \geq m_{k_s} \geq k_s \geq s$, it follows from 
Lemma \ref{lem-Omega} that $\mu_{y^{(s)}} \rightarrow \mu_{y^*}$ 
or equivalently $\mu_{y^{(s)}} \rightarrow \mu_y$ as 
$s \rightarrow \infty$. To see that
$\Phi (y) \leq \liminf_{s \rightarrow \infty} 
\Phi (y^{(s)})$, fix $s \geq 1$. By Lemma \ref{lem-Gamma},
\begin{equation*}
\begin{split}
\Phi (y^{(s)}) & = \lim_{p \rightarrow \infty} 
\frac{1}{p \cdot l_s} \ln S (y^{(s)}_0 \cdots y^{(s)}_{p \cdot l_s})
\geq \limsup_{p \rightarrow \infty} \frac{1}{p \cdot l_s} 
\ln \big[ \Gamma^{p \cdot l_s}_{y^{(s)}} (b_*, b_*) \big] \\
& \geq \limsup_{p \rightarrow \infty} \frac{1}{l_s} 
\ln \big[ \Gamma^{l_s}_{y^{(s)}} (b_*, b_*) \big] 
= \frac{1}{l_s} \ln \big[ \Gamma^{l_s}_{y^{(s)}} (b_*, b_*) \big]. 
\end{split}
\end{equation*}
It is clear that $\Gamma^{n_{k_s}}_{y^*} (b_*, c_*) 
\leq \Gamma^{l_s}_{y^{(s)}} (b_*, b_*)$. 
Thus from (\ref{eq-Phi=lim}), 
\begin{equation*}
\begin{split}
\Phi (y) & \leq \liminf_{s \rightarrow \infty} \frac{1}{n_{k_s}}  
\ln \big[ \Gamma^{n_{k_s}}_{y^*} (b_{k_s}, c_{k_s}) \big] 
= \liminf_{s \rightarrow \infty} \frac{1}{n_{k_s}} 
\ln \big[ \Gamma^{n_{k_s}}_{y^*} (b_*, c_*) \big] \\ 
& \leq \liminf_{s \rightarrow \infty} \frac{1}{l_s} 
\ln \big[ \Gamma^{l_s}_{y^{(s)}} (b_*, b_*) \big] 
\leq \liminf_{s \rightarrow \infty} \Phi (y^{(s)}),
\end{split}
\end{equation*}
which completes the proof.
\end{proof}

Let $E = \{ y \in \cg(Y) \vert \int \ct d \mu_y = \ct (y) \}$.
Then $\nu(E)=1$ for every ergodic invariant measure $\nu$ on $Y$, and
hence $\nu(E)=1$ for every $\nu \in M(Y)$. 
For $y \in E$, let $\{ y^{(s)} \}^\infty_{s = 1} \subset P (Y)$ be a
sequence obtained from Lemma \ref{lem-mu_Phi} so that
$\mu_{y^{(s)}} \rightarrow \mu_y$ as $s \rightarrow \infty$
and $\Phi (y) \leq \liminf_{s \rightarrow \infty} \Phi (y^{(s)})$.
By Lemma \ref{lem-cl_usc},
\begin{equation*}
\limsup_{s \rightarrow \infty} \ct (y^{(s)}) 
= \limsup_{s \rightarrow \infty} \int\!\ct d \mu_{y^{(s)}} 
\leq \int\!\ct d \mu_y = \ct (y).
\end{equation*}
It follows from Lemma \ref{lem-P(Y)} that
$\ct (y^{(s)}) = \Phi (y^{(s)})$ for all $s \geq 1$. Thus
\begin{equation*}
\Phi (y) \leq \liminf_{s \rightarrow \infty} \Phi (y^{(s)}) 
= \liminf_{s \rightarrow \infty} \ct (y^{(s)}) 
\leq \ct (y) \leq \Phi (y).
\end{equation*}
Hence $\ct (y) = \Phi (y)$ for all $y \in E$. Let $\nu \in M (Y)$.
Since $\nu (E) = 1$, we have
\begin{equation*}
\int\!\ct \h d \nu = \int_E\!\ct \h d \nu = \int_E\!\Phi \h d \nu
= \int\!\Phi \h d \nu ,
\end{equation*} 
which completes the proof of Theorem \ref{thm-main}.

To extend Theorem \ref{thm-main} to the case of an arbitrary potential
$f \in C(X)$, 
we need to define a suitable function corresponding to $\Phi_f$.
Let $f \in C (X)$. Fix $n \geq 1$. For each $i = 0, 1, \cdots, n-1$,
define a block map $F_n^i : \cb_n (X) \rightarrow \Br$ by
\begin{equation*} 
F_n^i (b_1 \cdots b_n) = \inf_{\sigma^{-i} x \in [b_1 \cdots b_n]} 
\exp (f (x))
\end{equation*}
(so the infimum is taken over all $x$ in the cylinder set
$_{-i}[b_1 \cdots b_n]_{n-i-1}$). 
For each $n$-block $B \in \cb_n (X)$, put  
\begin{equation*}
s_f (B) = \prod^{n-1}_{i=0} F_n^i (B),
\end{equation*}
and for a block $C \in \cb_n (Y)$ put 
\begin{equation*}
S_f (C) = \sum_{\pi(B)=C} s_f (B).
\end{equation*} 
Define $\Psi_f : Y \rightarrow \Br$ by
\begin{equation*}
\Psi_f (y) = \limsup_{n \rightarrow \infty} \frac{1}{n} 
\ln S_f (y_0 y_1 \cdots y_{n-1}) \hspace{0.3in} \text{for } y \in Y. 
\end{equation*}

Now, given $n \geq 1$, for each $i = 0, 1, \cdots, n-1$
define $\widetilde{F}_n^i : \cb_n (X) \rightarrow \Br$ by
\begin{equation*} 
\widetilde{F}_n^i (b_1 \cdots b_n) 
= \sup_{\sigma^{-i} x \in [b_1 \cdots b_n]} \exp (f (x)).
\end{equation*}
Using $\widetilde{F}_n^i$, 
we similarly define $\widetilde{s}_f, \widetilde{S}_f$ and 
$\widetilde{\Psi}_f$ as follows. For each $B \in \cb_n (X)$,
\begin{equation*}
\widetilde{s}_f (B) = \prod^{n-1}_{i=0} \widetilde{F}_n^i (B),
\end{equation*}
and for $C \in \cb_n (Y)$,
\begin{equation*}
\widetilde{S}_f (C) = \sum_{\pi(B)=C} \widetilde{s}_f (B).
\end{equation*} 
Define $\widetilde{\Psi}_f : Y \rightarrow \Br$ by
\begin{equation*}
\widetilde{\Psi}_f (y) = \limsup_{n \rightarrow \infty} \frac{1}{n} 
\ln \widetilde{S}_f (y_0 y_1 \cdots y_{n - 1}) \hspace{0.3in} 
\text{for } y \in Y.
\end{equation*}

It can be easily shown that  
if $f$ depends on just $2$ coordinates, $x_0x_1$, of each point 
$x \in X$, then $\Psi_f$ and $\widetilde{\Psi}_f$ both are equivalent
to $\Phi_f$,
 and hence Theorem \ref{thm-main} would apply. 
In 
the current more general situation, we still have the analogous result.

\begin{thm}
\label{thm_ext}
Let $X$ be an irreducible shift of finite type, $Y$ a subshift, 
$\pi : X \rightarrow Y$ a factor map, and $f \in C (X)$. 
For each $\nu \in M (Y)$, 
we have $P(\pi,f)(y) = \Psi_f(y) = \widetilde{\Psi}_f(y)$ a.e. 
$d\nu (y)$, and hence
\begin{equation*}
\int\!P (\pi, f) \h d \nu = \int\!\Psi_f \h d \nu 
= \int\!\widetilde{\Psi}_f \h d \nu.
\end{equation*}
\end{thm}

\begin{proof}
As before, we may assume that $f \geq 0$. 
Thus we have a constant $M > 0$ such that $1 \leq \exp (f) \leq M$
and hence $1 \leq F_n^i \leq \widetilde{F}_n^i \leq M$ 
for all $n \geq 1$ and for $i = 0, 1, \cdots, n-1$. 
Let $\ct = P (\pi, f)$ as before.
Set $s = s_f$, $S = S_f$, and $\Psi = \Psi_f$. For $y \in Y$, define
\begin{equation*}
\theta (y) = \limsup_{n \rightarrow \infty} \frac{1}{n} \ln 
\bigg[ \sum_{x \in D_n (y)} s (x_0 x_1 \cdots x_{n - 1}) \bigg].
\end{equation*} 
We show that $\ct (y) = \theta (y)$ for all $y \in Y$. 
Note first that if $x_0 \cdots x_{n-1} \in \cb_n (X)$, 
$n \geq 1$, then 
\begin{equation*}
\begin{split}
s (x_0 \cdots x_{n-1}) & = \prod^{n-1}_{i=0} 
\inf_{\sigma^{-i} z \in [x_0 \cdots x_{n-1}]} \exp (f (z)) \\
& = \prod^{n-1}_{i=0} \inf_{z \in [x_0 \cdots x_{n-1}]} 
\exp (f (\sigma^i z)).  
\end{split}
\end{equation*}
Thus from (\ref{eq_P(pi,f)}), given $y \in Y$,
\begin{equation*}
\begin{split}
\theta (y) & = \limsup_{n \rightarrow \infty} \frac{1}{n} \ln \bigg[ 
\sum_{x \in D_n (y)} \prod^{n-1}_{i=0} 
\inf_{z \in [x_0 \cdots x_{n-1}]} \exp (f (\sigma^i z)) \bigg] \\
& \leq \limsup_{n \rightarrow \infty} \frac{1}{n} 
\ln \bigg[ \sum_{x \in D_n (y)} \prod^{n-1}_{i=0}
\exp \big( f (\sigma^i x) \big) \bigg] = \ct (y).
\end{split}
\end{equation*} 

Suppose that there exist $y \in Y$ and $\epsilon > 0$ for which
\begin{equation*}
\theta (y) = \ct (y) - 2 \epsilon.
\end{equation*} 
Since $f$ is uniformly continuous, there is $p \geq 0$ such that
whenever $\rho (z, x) < 2^{-p}$ or,
 equivalently, 
$z_{[-p, p]} = x_{[-p, p]}$, then  
$\vert f (z) - f (x) \vert < \epsilon$.
Fix $x \in X$ and $n > 2p$. For each $i = p, p+1, \cdots, n-p-1$, 
\begin{equation*} 
\begin{split}
F_n^i (x_0 x_1 \cdots x_{n-1}) 
& = \inf_{z \in [x_0 \cdots x_{n-1}]} \exp (f (\sigma^i z)) \\
& \geq \inf_{\substack{z \in X \\ 
\rho (\sigma^i z, \sigma^i x) < 2^{-p}}} \exp (f (\sigma^i z)) 
\geq \exp \big( f (\sigma^i x) - \epsilon \big). 
\end{split} 
\end{equation*}
Since $F^i_n \geq 1$ for each $i$ and $\exp (f) \leq M$, we have 
\begin{equation*} 
\begin{split}
s (x_0 x_1 \cdots x_{n-1}) 
& \geq \prod^{n-p-1}_{i=p} F_n^i (x_0 x_1 \cdots x_{n-1}) 
\geq \prod^{n-p-1}_{i=p} \exp \big( f (\sigma^i x) - \epsilon \big) \\ 
& \geq M^{-2p} e^{-(n-2p) \epsilon}
\prod^{n-1}_{i=0} \exp (f (\sigma^i x)). 
\end{split} 
\end{equation*}
It follows that
\begin{equation*}
\begin{split}
\theta (y) & \geq \limsup_{n \rightarrow \infty} \frac{1}{n} \ln 
\bigg[ M^{-2p} e^{-(n-2p) \epsilon} \sum_{x \in D_n (y)} 
\prod^{n-1}_{i=0} \exp (f (\sigma^i x)) \bigg] \\
& = \limsup_{n \rightarrow \infty} \frac{1}{n} \ln \bigg[ 
\sum_{x \in D_n (y)} \prod^{n-1}_{i=0} \exp (f (\sigma^i x)) \bigg] 
- \epsilon = \ct (y) - \epsilon,
\end{split}
\end{equation*} 
which is a contradiction. Therefore $\ct (y) = \theta (y)$.

Observe next that for each $y \in Y$,
\begin{equation*}
\begin{split}
\theta (y) & = \limsup_{n \rightarrow \infty} \frac{1}{n} \ln \bigg[ 
\sum_{x \in D_n (y)} s (x_0 x_1 \cdots x_{n - 1}) \bigg] \\
& \leq \limsup_{n \rightarrow \infty} \frac{1}{n} \ln \bigg[ 
\sum_{\substack{\pi (x_0 \cdots x_{n-1}) \\ = y_0 \cdots y_{n-1}}} 
s (x_0 x_1 \cdots x_{n-1}) \bigg] = \Psi (y) ,
\end{split} 
\end{equation*}
since the summation for $\Psi (y)$ is taken over the larger set
than the one for $\theta (y)$. 
As we have seen before, it is possible to have 
$\theta (y) < \Psi (y)$ for some $y \in Y$. We will prove
that if $y \in Y$ is periodic, then $\theta (y) = \Psi (y)$, or, 
equivalently $\ct (y) = \Psi (y)$. 

Fix $y \in Y$ and for each $n \geq 1$, let
\begin{equation*}
\tau_n (y) = \sum_{x \in D_n(y)} s (x_0 x_1 \cdots x_{n-1}).  
\end{equation*}
Assuming that $D_n (y) \subset D_{n+1} (y)$ for each $n \geq 1$
(without loss of generality), we have
\begin{equation*} 
\begin{split}
\tau_n (y) & = \sum_{x \in D_n (y)} \prod^{n-1}_{i=0} 
\inf_{z \in x_{[0, n-1]}} \exp (f (\sigma^i z)) 
\leq \sum_{x \in D_n (y)} \prod^{n-1}_{i=0} 
\inf_{z \in x_{[0, n]}} \exp (f (\sigma^i z)) \\
& \leq \sum_{x \in D_{n+1} (y)} \prod^n_{i=0} \inf_{z \in x_{[0, n]}}  
\exp (f (\sigma^i z)) = \tau_{n+1} (y).
\end{split} 
\end{equation*}
Hence $\tau_n (y)$ increases as $n$ increases.
An easy computation shows that given an increasing sequence 
$\{ a_n \}$, for any $q \geq 1$,
$\limsup_{n \rightarrow \infty} a_n / n 
= \limsup_{n \rightarrow \infty} a_{nq} / (nq)$. Letting
$a_n = \ln \tau_n (y)$ proves that for any $q \geq 1$,
\begin{equation}
\label{eq_theta(y)}
\theta (y) = \limsup_{n \rightarrow \infty} \frac{1}{nq} \ln 
\bigg[ \sum_{x \in D_{nq} (y)} s (x_0 x_1 \cdots x_{nq-1}) \bigg].
\end{equation}

Let $y \in P_q (Y)$, $q \geq 1$. We claim that
\begin{equation}
\label{eq_Psi(y)}
\begin{split}
\Psi (y) & = \limsup_{n \rightarrow \infty} \frac{1}{nq} 
\ln S (y_0 y_1 \cdots y_{nq-1}) \\
& = \limsup_{n \rightarrow \infty} \frac{1}{nq} \ln \bigg[ 
\sum_{\substack{\pi (x_0 \cdots x_{nq-1}) \\ = y_0 \cdots y_{nq-1}}} 
s (x_0 x_1 \cdots x_{nq-1}) \bigg].
\end{split}
\end{equation}
Then, using (\ref{eq_theta(y)}) and (\ref{eq_Psi(y)}),
we can proceed as in Lemma \ref{lem-P(Y)} to show that 
$\Psi (y) \leq \theta (y)$ and therefore $\Psi (y) = \ct (y)$. 

To verify (\ref{eq_Psi(y)}), let $\epsilon > 0$ be given and  
(using uniform continuity of $f$) choose $p \geq 0$ so that whenever 
$\rho (z, x) < 2^{-p}$ or $z_{[-p, p]} = x_{[-p, p]}$, then  
$\vert f (z) - f (x) \vert < \epsilon$.
Given an integer $m \geq 2p + q$, there is $n \in \Bn$ such that
$2p < n q \leq m < (n+1) q$. Put $k = m - nq \geq 0$ and
$v = y_0 y_1 \cdots y_{nq-1}$. Then 
\begin{equation} 
\begin{split}
\label{eq_S(v)}
S (y_0 y_1 \dots y_{m-1}) & = \sum_{\pi (u) = y_0 \cdots y_{m-1}} 
\prod^{m-1}_{i=0} \inf_{z \in [u]} \exp (f (\sigma^i z)) \\
& \leq \sum_{\pi (u) = v} 
\sum_{\substack{w \in \cb_k (X) \\ u w \in \cb (X)}} 
\prod^{m-1}_{i=0} \inf_{z \in [u w]} \exp (f (\sigma^i z)) \\
& \leq \sum_{\pi (u) = v} 
\sum_{\substack{w \in \cb_k (X) \\ u w \in \cb (X)}} 
M^{k+2p} \prod^{nq-p-1}_{i=p} \inf_{z \in [uw]} \exp (f(\sigma^i z)). 
\end{split}
\end{equation}
Let $u \in \cb_{nq} (X)$ and $w \in \cb_k (X)$ be any blocks such
that $u w \in \cb (X)$. Let $p \leq i \leq nq-p-1$. Note that if
$z \in [u w]$ and $\bar{z} \in [u]$, then
$\rho (\sigma^i z, \sigma^i \bar{z}) < 2^{-p}$ so that
$\vert f (\sigma^i z) - f (\sigma^i \bar{z}) \vert < \epsilon$. Thus
\begin{equation*}
\inf_{z \in [u w]} \exp (f (\sigma^i z)) 
\leq \inf_{z \in [u]} \exp \big( f (\sigma^i z) + \epsilon \big).  
\end{equation*}
From this inequality and (\ref{eq_S(v)}), we have
\begin{equation*}
\begin{split}
S (y_0 y_1 \dots y_{m-1}) & \leq M^{k+2p} \sum_{\pi (u) = v} 
\sum_{\substack{w \in \cb_k (X) \\ u w \in \cb (X)}} 
\prod^{nq-p-1}_{i=p} \inf_{z \in [u]} 
\exp \big( f (\sigma^i z) + \epsilon \big) \\
& \leq M^{k+2p} \cdot \big\vert \cb_k (X) \big\vert 
\sum_{\pi (u) = v} \prod^{nq-1}_{i=0} \inf_{z \in [u]} 
\exp \big( f (\sigma^i z) + \epsilon \big) \\
& = M^{k+2p} \cdot \big\vert \cb_k (X) \big\vert \cdot 
e^{nq \epsilon} \cdot S (y_0 y_1 \dots y_{nq-1}).
\end{split} 
\end{equation*}
Since $k < q$ and $p, q$ are fixed, it follows that
\begin{equation*}
\begin{split}
\Psi (y) & = \limsup_{m \rightarrow \infty} \frac{1}{m} 
\ln S (y_0 \cdots y_{m-1}) \\
& \leq \limsup_{n \rightarrow \infty} \frac{1}{nq} 
\ln S (y_0 y_1 \dots y_{nq-1}) + \epsilon.
\end{split}
\end{equation*}
Taking $\epsilon$ arbitrarily small completes the proof of the
claim. 

Next, we will show that $\widetilde{\Psi}_f (y) = \ct (y)$ for 
$y \in P (Y)$. It is straightforward to check that 
$\ct (y) \leq \widetilde{\Psi}_f (y)$ for all $y \in Y$.
Set $\widetilde{s} = \widetilde{s}_f$, 
$\widetilde{S} = \widetilde{S}_f$, and 
$\widetilde{\Psi} = \widetilde{\Psi}_f$. 
Define $\widetilde{\theta} : Y \rightarrow \Br$ by  
\begin{equation*}
\widetilde{\theta} (y) = \limsup_{n \rightarrow \infty} \frac{1}{n} 
\ln \bigg[ \sum_{x \in D_n (y)} 
\widetilde{s} (x_0 x_1 \cdots x_{n - 1}) \bigg].
\end{equation*} 
Note that $\ct (y) \leq \widetilde{\theta} (y)$. Suppose there
exist $y \in Y$ and $\epsilon > 0$ for which
\begin{equation*}
\widetilde{\theta} (y) = \ct (y) + 2 \epsilon.
\end{equation*} 
Similarly
 to the foregoing, there is $p \geq 0$ such that
whenever $\rho (z, x) < 2^{-p}$ or equivalently, 
$z_{[-p, p]} = x_{[-p, p]}$, then  
$\vert f (z) - f (x) \vert < \epsilon$.
Fix $x \in X$ and $n > 2p$. If $p \leq r \leq n-p-1$, then
\begin{equation*} 
\begin{split}
\widetilde{F}_n^r (x_0 x_1 \cdots x_{n-1}) 
& = \sup_{\sigma^{-r} z \in [x_0 \cdots x_{n-1}]} \exp (f (z)) \\
& \leq \sup_{\rho (z, \sigma^r x) < 2^{-p}} \exp (f (z)) 
\leq \exp (f (\sigma^r x) + \epsilon).
\end{split} 
\end{equation*}
Thus
\begin{equation*} 
\begin{split}
\widetilde{s} (x_0 x_1 & \cdots x_{n-1}) \\ 
& = \prod^{p-1}_{r=0} \widetilde{F}_n^r (x_0 \cdots x_{n-1}) 
\prod^{n-p-1}_{r=p} \widetilde{F}_n^r (x_0 \cdots x_{n-1}) 
\prod^{n-1}_{r=n-p} \widetilde{F}_n^r (x_0 \cdots x_{n-1}) \\ 
& \leq M^{2p} \prod^{n-p-1}_{r=p} 
\widetilde{F}_n^r (x_0 \cdots x_{n-1})  
\leq M^{2p} \prod^{n-p-1}_{r=p} \exp (f (\sigma^r x) + \epsilon) \\
& \leq M^{2p} e^{(n-2p) \epsilon}
\exp \Big( \sum^{n-p-1}_{r=p} f (\sigma^r x) \Big) 
\leq M^{2p} e^{(n-2p) \epsilon}
\exp \Big( \sum^{n-1}_{r=0} f (\sigma^r x) \Big). 
\end{split} 
\end{equation*}
It follows that
\begin{equation*}
\begin{split}
\widetilde{\theta} (y) 
& \leq \limsup_{n \rightarrow \infty} \frac{1}{n} \ln 
\bigg[ M^{2p} e^{(n-2p) \epsilon} \sum_{x \in D_n (y)} 
\exp \Big( \sum^{n-1}_{r=0} f (\sigma^r x) \Big) \bigg] \\
& = \epsilon + \limsup_{n \rightarrow \infty} \frac{1}{n} \ln 
\bigg[ \sum_{x \in D_n (y)} 
\exp \Big( \sum^{n-1}_{r=0} f (\sigma^r x) \Big) \bigg] 
= \epsilon + \ct (y) ,
\end{split}
\end{equation*} 
which is a contradiction. Therefore 
$\ct (y) = \widetilde{\theta} (y)$.

Let $y \in P_q (Y)$, $q \geq 1$. It is not difficult to see that 
\begin{equation*}
\ct (y) = \limsup_{n \rightarrow \infty} \frac{1}{nq} \ln \bigg[ 
\sum_{x \in D_{nq} (y)} \widetilde{s} (x_0 x_1 \cdots x_{nq - 1})
\bigg] = \widetilde{\theta} (y)
\end{equation*}
and 
\begin{equation*}
\widetilde{\Psi} (y) = \limsup_{n \rightarrow \infty} \frac{1}{nq} 
\ln \widetilde{S} (y_0 \cdots y_{n q - 1}).
\end{equation*}
We can proceed again as in Lemma \ref{lem-P(Y)} to show that 
$\widetilde{\Psi} (y) = \ct (y)$. 
The remainder of the proof is the same as before.
\end{proof}

Recall that according to Theorem 4.6 of \cite{Wal}, if for each
$n = 1, 2, \cdots$ and $y \in Y$ we denote by $D_n (y)$ a set
consisting of exactly one point from each nonempty set 
$[x_0 \cdots x_{n-1}] \cap \pi^{-1} (y)$, then 
for each $f \in C(Y)$,
\begin{equation*}
P(\pi, f)(y) = \limsup_{n \rightarrow \infty} \frac{1}{n} 
\ln \bigg[ \sum_{x \in D_n (y)} \exp \Big( \sum_{i=0}^{n-1} 
f (\sigma^i x) \Big) \bigg].
\end{equation*}
From Theorem \ref{thm_ext} it now follows that we will obtain the
value of $P(\pi,f) (y)$ a.e. with respect to every invariant measure
on $Y$ if we delete from the definition of $D_n (y)$ the requirement
that $x \in \pi^{-1} (y)$:

\begin{cor*}
For each $n = 1, 2, \cdots$ and $y \in Y$ denote by $E_n (y)$ a set
consisting of
exactly one point from each nonempty cylinder 
$[x_0 \cdots x_{n-1}] \subset \pi^{-1} [y_0 \cdots y_{n-1}]$. 
Then for each $f \in C(Y)$, 
\begin{equation*}
P(\pi,f) (y) = \limsup_{n \rightarrow \infty} \frac{1}{n} 
\ln \bigg[ \sum_{x \in E_n (y)} \exp \Big( \sum_{i=0}^{n-1} 
f (\sigma^i x) \Big) \bigg]
\end{equation*}
a.e. with respect to every invariant measure on Y.
\end{cor*}

\begin{ack*}
The authors thank the University of Warwick, where much of this
research was accomplished, for its outstanding hospitality.
\end{ack*}

\bibliographystyle{amsplain}

\end{document}